# Formule d'Itô pour des diffusions uniformément elliptiques, et processus de Dirichlet

K.Dupoiron[*]/P.Mathieu[†]/J.San Martin[‡]

15 novembre 2018

**Mots clés :** Processus de Dirichlet, formes de Dirichlet, formule d'Itô, intégrales stochastiques, covariations quadratiques.

**Résumé.** - Soit $X$ une diffusion uniformément elliptique sur $\mathbb{R}^d$, $F$ une fonction dans $H^1_{\text{loc}}(\mathbb{R}^d)$ et $\nu$ la loi initiale de la diffusion. On montre que si l'intégrale $\int |\nabla F|^2(x) U\nu(x)\, dx$ est finie, où $U\nu$ désigne le potentiel de la mesure $\nu$, alors $F(X)$ est un processus de Dirichlet. Si de plus, $F$ appartient à $H^2_{\text{loc}}(\mathbb{R}^d)$ et si les intégrales $\int |\nabla F|^2(x) U\nu(x)\, dx$ et $\int |\nabla f_k|^2(x) U\nu(x)\, dx$ sont finies, pour les dérivées faibles $f_k$ de $F$, alors on peut écrire une formule d'Itô. En particulier, on définit l'intégrale progressive $\int \nabla F(X)\, dX$ et les covariations quadratiques $[f_k(X), X^k]$.

## 1 Introduction

Soit $(X_t)_{t\geq 0}$ un processus de diffusion à valeurs dans $\mathbb{R}^d$, soit $F : \mathbb{R}^d \to \mathbb{R}$ une fonction borélienne et $Y := F(X)$. Si $F$ est de classe $C^2$, on sait décrire le comportement du processus $Y$ grâce à la formule d'Itô, pourvu que les coefficients du générateur de la diffusion soient suffisamment réguliers. Dans ce cas $Y$ est une semi-martingale. Cependant si $F$ ou si les coefficients de la diffusion sont moins réguliers, le processus $Y$ n'est plus en général une semi-martingale. Que peut-on dire de $Y$ ? A-t-on encore une formule d'Itô ? Pour répondre à ces questions, il faut recourir à la notion de processus de Dirichlet, introduite par Fukushima ([2]), qui généralise celle de semi-martingale.

**Définition 1 : Processus de Dirichlet**
Soit $S$ un réel strictement positif, et $(A_t,\ 0 \leq t \leq S)$ un processus stochastique.

---

[*]CMI, Université de Provence, 39 rue Joliot Curie, 13453 Marseille Cedex 13, France, Karelle.Dupoiron@cmi.univ-mrs.fr

[†]CMI, Université de Provence, 39 rue Joliot Curie, 13453 Marseille Cedex 13, France, Pierre.Mathieu@cmi.univ-mrs.fr

[‡]Universidad de Chile, Facultad de Ciencias Físicas y Matemáticas, Departamento de Ingeniería Matemática, Casilla 170-3 Correo 3, Santiago, Chile, jsanmart@dim.uchile.cl

On dit que $A$ est un processus de Dirichlet, si on peut trouver une décomposition $A_t = M_t + N_t$ de $A$, telle que
-$M$ est une martingale continue de carré intégrable,
-$N$ est un processus à variations quadratiques nulles, c'est-à-dire :

$$\lim_{n \to +\infty} \sum_{t_i \in D_n} \left(N_{t_{i+1}} - N_{t_i}\right)^2 = 0$$

en probabilité, où $D_n$ désigne la subdivision dyadique d'ordre $n$ de $[0, S]$.

Supposons donnée une forme de Dirichlet symétrique $\mathcal{E}$ de domaine $\mathcal{D}(\mathcal{E})$ inclus dans $L^2(\mathbb{R}^d, dx)$, et $(X_t, P_t)$ le processus de Hunt associé. On note $\mathbb{P}_x$ la loi de $X$ lorsque $X_0 = x$, $\mathbb{P}_\nu$ la loi de $X$ lorsque la loi de $X_0$ est $\nu$. Alors Fukushima [2] montre que si $u \in \mathcal{D}(\mathcal{E})$, et si $\tilde{u}$ est une version quasi-continue de $u$, la fonctionnelle additive $A_t^{(u)} := \tilde{u}(X_t) - \tilde{u}(X_0)$ est un processus de Dirichlet, lorsque la loi $\mu$ de $X_0$ est absolument continue par rapport à la mesure de Lebesgue. Notons alors $A^{(u)} = M^{(u)} + N^{(u)}$ sa décomposition en tant que processus de Dirichlet. La partie martingale $M^{(u)}$ est une martingale sous $P_x$ pour quasi tout $x$, et $N^{(u)}$ est d'énergie nulle au sens où :
$\lim_{t \to 0} \dfrac{1}{2t} E_m(N_t^2) = 0$ où $m$ désigne la mesure de Lebesgue. En particulier $N^{(u)}$ est à variations quadratiques nulles au sens décrit ci-dessus par la propriété de Markov, sous la mesure $\mathbb{P}_m$. Si de plus la fonction $u$ est dans le domaine $L^2$ du générateur de $X$, alors $N^{(u)}$ est à variations bornées et $A^{(u)}$ est une semi-martingale.

**Cadre d'étude :** Soit $a(x) = (a_{ij}(x))_{1 \leq i,j \leq d}$ une famille de fonctions mesurables, localement intégrables sur $\mathbb{R}^d$, telles que $a_{ij} = a_{ji}$ pour tous $1 \leq i, j \leq d$. On note $\sigma$ la racine carrée de $a$. On suppose de plus qu'il existe une constante $\lambda \geq 1$ telle que

$$\frac{1}{\lambda} \sum_{i=1}^d \xi_i^2 \leq \sum_{i,j=1}^d a_{ij}(x)\xi_i\xi_j \leq \lambda \sum_{i=1}^d \xi_i^2, \qquad \forall \xi \in \mathbb{R}^d, \forall x \in \Omega.$$

On note $L$ l'opérateur $\mathrm{div}(a \,\mathrm{grad})$, et $\mathcal{E}$ la forme de Dirichlet associée : $(\mathcal{E}, \mathcal{D}(\mathcal{E}))$ est la fermeture dans $L^2(\mathbb{R}^d, dx)$ de la forme bilinéaire, fermable

$$\mathcal{E}(f, g) = \int_{\mathbb{R}^d} -Lf(x)g(x)dx = \int_{\mathbb{R}^d} \sigma\nabla f(x).\sigma\nabla g(x)dx$$

définie pour toutes fonctions $f$ et $g \in C_c^\infty$. Sous les hypothèses faites sur $\sigma$, $\mathcal{D}(\mathcal{E}) = H^1(\mathbb{R}^d)$. Soit $(X_t, P_t)$ le processus de Hunt associé.

**Notations 1** Pour une loi $\nu$ donnée, on écrira :

$\nu_t(x) := \dfrac{dP_t\nu}{dx}(x)$  $\quad$ $P_t\nu$ étant absolument continue par rapport à la mesure de Lebesgue, pour $t > 0$.

$U\nu(x) := \displaystyle\int_0^{+\infty} ds \, e^{-s}\nu_s(x)$  $\quad$ lorsque cette intégrale converge.

Föllmer et Protter [1] ont démontré que dans le cas où $X$ est un mouvement brownien ($a = Id$), et lorsque $F \in H^1_{\text{loc}}$, en dimension $d$, on peut obtenir une formule d'Itô généralisée où le terme d'ordre 2 est remplacé par les covariations quadratiques $[f_k(X), X]$ où $f_k = \frac{\partial F}{\partial x_k}$. Celles-ci sont données comme limites en probabilité de sommes de Riemann, pour tout point de départ en dehors d'un ensemble de capacité nulle.

On montre ici, par un calcul direct, que si $F$ vérifie la condition d'intégrabilité

$$\int_{\mathbb{R}^d} dz \, |\nabla F(z)|^2 U\nu(z) < +\infty \tag{1}$$

étant donnée une loi initiale $\nu$, alors $Y = F(X)$ est un processus de Dirichlet. Si de plus les dérivées de $F$ vérifient la condition d'intégrabilité

$$\int_{\mathbb{R}^d} dz \, |\nabla \partial_k F(z)|^2 U\nu(z) < +\infty \tag{2}$$

on peut alors écrire une formule d'Itô sous $\mathbb{P}_\nu$ :

$$F(X_t) - F(X_0) = \int_0^t \nabla F(X_s) dX_s + \frac{1}{2} \sum_{k=1}^d [f_k(X), X]_t \quad \mathbb{P}_\nu - p.s. \tag{3}$$

où

$$\int_0^t \nabla F(X_s) dX_s := \lim_{n \to +\infty} \sum_{t_i \in D_n} \nabla F(X_{t_i}).(X_{t_{i+1}} - X_{t_i}) \tag{4}$$

$$\text{et } [f_k(X), X]_t := \lim_{n \to +\infty} \sum_{t_i \in D_n} (f_k(X_{t_{i+1}}) - f_k(X_{t_i}))(X^k_{t_{i+1}} - X^k_{t_i}) \tag{5}$$

$D_n$ désignant la suite des subdivisions dyadiques de $[0, t]$, les limites étant prises en probabilité, sous $\mathbb{P}_\nu$. On en déduit une formule d'Itô, sous $\mathbb{P}_{x_0}$ pour quasi tout point de départ $x_0$, lorsque $F$ est dans $H^2_{\text{loc}}$, les limites (4) et (5) existant en probabilité sous $\mathbb{P}_{x_0}$. Comme dans le cas brownien, elles sont obtenues pour quasi tout point de départ $x_0$ et non pour tout point de départ. On insiste cependant sur deux points.

D'une part l'égalité (4) est cette fois une définition : $X$ n'est pas une martingale, l'intégrale (4) n'est donc pas une intégrale stochastique. En particulier la décomposition (3) ne correspond pas à la décomposition donnée par Fukushima [2] en somme de martingale et de processus à variations quadratiques nulles, contrairement au cas brownien. Toutefois on peut encore dire que la covariation $[f_k(X), X]_t$ est à variations quadratiques nulles. L'intégrale progressive définie par (4) contient donc une partie martingale et une partie à variations quadratiques nulles.

D'autre part la condition d'intégrabilité (1) est la condition minimale requise pour que $F(X)$ soit un processus de Dirichlet, partant d'une loi $\nu$, puisque le crochet de la partie martingale est alors $\int_0^S |\sigma \nabla F|^2(X_s)ds$. Comme dans [5], on note que c'est le lien entre la régularité de $F$ et la loi initiale $\nu$ qui donne la propriété de Dirichlet. Cette dépendance en $\nu$ disparaît, si on demande une intégrabilité plus forte sur $F$. En particulier si $F \in W^{1,p}$ pour un $p > d \vee 2$, alors la condition (1) est vérifiée pour toute loi $\nu$ (appliquer l'inégalité de Hölder), et on retrouve ainsi un résultat de Rozkosz [6] : $Y = F(X)$ est alors un processus de Dirichlet pour toute loi initiale.

Enfin, on notera que la seule condition (1) est requise pour $F$ dans la première partie (pour ses dérivées également dans la seconde partie) : elle correspond à la norme $\|F\|_2(x_0)$ de Föllmer et Protter [1] lorsque $\nu$ est la mesure de Dirac $\delta_{x_0}$ en $x_0$. En revanche il n'est pas nécessaire d'imposer une condition d'intégrabilité supplémentaire de type $\|F\|_1(x_0)$. Celle-ci apparaît dans [1] au moment de prouver la convergence des sommes rétrogrades $\sum F(X_{t_{i+1}})(X_{t_{i+1}} - X_{t_i})$. La preuve de Föllmer et Protter utilise une formule explicite pour le mouvement brownien retourné dans le temps, qui n'existe pas en général. Ici, on prouve directement la convergence des sommes (5), sans démontrer, séparément, la convergence des sommes progressives et des sommes rétrogrades.

## 2   Diffusion elliptique et processus de Dirichlet

### 2.1   Définitions et résultats

On se place dans le cadre décrit au paragraphe 1.

**Notations 2 :** Sauf mention contraire, $\nu$ désignera une mesure de probabilité, qui peut charger des ensembles de capacité nulle.
$(D_n)$ sera la suite des subdivisions dyadiques de l'intervalle $[0, S]$, $S$ étant un réel fixé.

**Proposition 1** *Soit $F$ une fonction mesurable dans $L^2(\mathbb{R}^d, dx)$, vérifiant (1). Alors sous $\mathbb{P}_\nu$, le processus $(Y_t := F(X_t), 0 \leq t \leq S)$ est un processus de Dirichlet. De plus on a pour $n$ assez grand :*

$$E_\nu \sum_{t_i \in D_n} \left(Y_{t_{i+1}} - Y_{t_i}\right)^2 \leq C \int_{\mathbb{R}^d} dz \, |\nabla F(z)|^2 U\nu(z)$$

*où $C$ est une constante qui dépend de $\lambda$, $d$ et $S$.*

Il faut d'abord remarquer que le processus $(Y_t)_{t>0}$ est "bien défini" pour les temps strictement positifs. En effet, si $F$ vérifie la condition (1), alors $F \in H^1_{\text{loc}}$. Ceci résulte de ce que pour tout compact $K$, $\inf_{x \in K} U\nu(x)$ est atteint et est strictement positif. Par conséquent, d'après [2], on peut supposer $F$ quasi-continue . Soient $F_1$ et $F_2$ deux versions quasi-continues de $F$. Pour tout $s > 0$, $P_s \nu$ est absolument continue par rapport à la mesure de Lebesgue, donc ne charge aucun ensemble de capacité nulle. On a alors
$$\mathbb{P}_\nu \left( \exists\, t > s,\ F_1(X_t) \neq F_2(X_t) \right) = 0.$$
grâce à la propriété de Markov. Ainsi $(F(X_t))_{t>0}$ est "bien défini".

En revanche, il faut pouvoir définir $F(X_0)$ sans ambiguïté : la mesure $\nu$ peut en effet charger des ensembles de capacité nulle, si bien que deux versions quasi-continues de $F$ peuvent différer en $X_0$. Comme $t \mapsto Y_t$ est continu sous $\mathbb{P}_\nu$ pour $t > 0$, il semble alors naturel de choisir une version $\tilde{F}$ de $F$, quasi-continue, de sorte que $\tilde{F}$ soit définie $\nu$ presque-sûrement et $\tilde{F}(X_0) = Y_0 = \lim_{t \to 0} F(X_t)$, la limite étant prise dans $L^2(\mathbb{P}_\nu)$. L'inégalité (12) prouve que cette limite existe, puisque la suite $(F(X_t))_{t>0}$ est de Cauchy en 0 dans $L^2(\mathbb{P}_\nu)$. Posons donc :
$$\tilde{F}(x) := \lim_{t \to 0} P_t F(x) = \lim_{t \to 0} E_x(F(X_t))$$

Cette limite existe pour quasi tout $x$, par forte continuité du semi-groupe, et aussi $\nu$ presque-sûrement d'après l'inégalité (12). Plus précisément, l'ensemble
$$A = \{x,\ E_x(F(X_t)) \text{ ne converge pas lorsque } t \text{ tend vers } 0\}$$
est de mesure $\nu$ nulle et de capacité nulle. De plus on a : $\tilde{F}(x) = F(x)$ quasi-partout. Donc $\tilde{F}$ est quasi-continue et $F(X_t) = \tilde{F}(X_t)$ $\mathbb{P}_\nu$-p.s. pour tout $t > 0$. $\tilde{F}(X_0)$ est bien définie et il reste à montrer qu'on a encore $\lim_{t \to 0} \tilde{F}(X_t) = \tilde{F}(X_0)$ dans $L^2(\mathbb{P}_\nu)$. Or pour $t > 0$ et $r > 0$, on a :
$$E_\nu \left( E_{X_t}(F(X_r)) - F(X_r) \right)^2 \tag{6}$$
$$= \int \nu(dx)\, E_x \left( E_{X_t}(F(X_r)) - F(X_r) \right)^2$$
$$= \int \nu(dx)\, E_x \left( F(X_{r+t}) - F(X_r) \right)^2 \tag{7}$$

par la propriété de Markov. D'après l'inégalité (12), (7) peut être rendu arbitrairement petit pourvu que $r$ et $t$ soient suffisamment petits. Faisant tendre $r$ vers 0 dans (6), on en déduit que pour tout $\varepsilon > 0$ et $t$ petit :
$$E_\nu \left( \tilde{F}(X_t) - \tilde{F}(X_0) \right)^2 \leq \varepsilon$$

Ainsi $\tilde{F}(X_t)$ tend vers $\tilde{F}(X_0)$ dans $L^2(\mathbb{P}_\nu)$.
Désormais on notera $F$, au lieu de $\tilde{F}$, la version quasi-continue ainsi choisie.

La proposition 1 nous permet de retrouver en particulier un résultat de [2] où $\nu$ est absolument continue ; dans ce cas, si $F$ appartient à $H^1 = \mathcal{D}(\mathcal{E})$, alors $(F(X_t), 0 \leq t \leq S)$ est un processus de Dirichlet. Ce résultat est en fait utilisé dans la preuve de la proposition.

Ici, $P_t\nu$ est absolument continue pour tout $\nu$ et tout $t > 0$. Donc $(Y_s, s \geq t)$ est un processus de Dirichlet pour tout $t > 0$. Démontrer la proposition 1, c'est donc contrôler les variations quadratiques de $(Y_s, s \leq t)$ pour tout $t$ petit. Pour cela nous avons besoin de contrôler les probabilités de transition de $X$.

**Lemme 1 : Estimées d'Aronson et Nash**
*Sous les hypothèses faites sur $a$, on a : $\exists M(\lambda, d)$*

$$\frac{1}{Mt^{d/2}} \exp-\left(\frac{M|y-x|^2}{t}\right) \leq p_t(x,y) \leq \frac{M}{t^{d/2}} \exp-\left(\frac{|y-x|^2}{Mt}\right)$$

Pour la démonstration on se référera à [7] ( théorèmes I.1.15 et I.2 ).

∎

**Notations 3 :** On pose
$$p_t^M(x,y) := \frac{M}{t^{d/2}} \exp-\left(\frac{|y-x|^2}{Mt}\right)$$
$$\nu_t^M(x) := \frac{dP_t^M \nu}{dx}(x)$$

Notons alors que $\int_0^S \nu_t^M(x)\, dt \leq C\, U\nu(x)$ pour une constante $C$.

Donnons immédiatement deux applications de la proposition 1 et de ces estimées :

**Corollaire 1** *Si $F \in W^{1,p}$ pour un $p > d \vee 2$, alors $F(X)$ est un processus de Dirichlet pour tout point de départ $x_0$.*

*Preuve :* Soit $x_0 \in \mathbb{R}^d$, et $\nu$ la mesure de Dirac en $x_0$, alors d'après les estimées gaussiennes, on a pour une certaine constante $M$ :

$$U\nu(x) \leq M \int_0^{+\infty} dt\, \frac{e^{-t-\frac{|x_0-x|^2}{Mt}}}{t^{d/2}}$$

Considérons pour $\alpha > 0$ et $A > 0$ les intégrales :

$$\int_0^{+\infty} dt\, \frac{e^{-t-\frac{\alpha}{t}}}{t^{d/2}} = \int_0^A dt\, \frac{e^{-t-\frac{\alpha}{t}}}{t^{d/2}} + \int_A^{+\infty} dt\, \frac{e^{-t-\frac{\alpha}{t}}}{t^{d/2}}$$

Le changement de variable $u = \frac{\alpha}{t}$ dans la première intégrale du second membre donne :

$$\begin{aligned}
\int_0^A dt\, \frac{e^{-t-\frac{\alpha}{t}}}{t^{d/2}} &\leq \int_0^A dt\, \frac{e^{-\frac{\alpha}{t}}}{t^{d/2}} \\
&\leq \frac{1}{\alpha^{d/2-1}} \int_{\frac{1}{A}}^{+\infty} du\, u^{d/2-2} e^{-u}
\end{aligned} \quad (8)$$

Par ailleurs, l'inégalité $e^{-x} \leq \frac{1}{1+x}$, $\forall x \geq 0$ donne dans la seconde intégrale :

$$\begin{aligned}\int_A^{+\infty} dt\, \frac{e^{-t-\frac{\alpha}{t}}}{t^{d/2}} &\leq \int_A^{+\infty} dt \frac{1}{1+\frac{\alpha}{t}} \frac{e^{-t}}{t^{d/2}} \\ &\leq \frac{1}{A+\alpha} \int_A^{+\infty} dt\, \frac{e^{-t}}{t^{d/2-1}} \end{aligned} \quad (9)$$

Soit $q > 1$ un entier, en utilisant (8) et (9), puis l'inégalité de Hölder :

$$\left(\int_0^{+\infty} dt\, \frac{e^{-t-\frac{\alpha}{t}}}{t^{d/2}}\right)^q \leq C\left(\left(\frac{1}{\alpha^{d/2-1}}\right)^q + \left(\frac{1}{A+\alpha}\right)^q\right)$$

pour une certaine constante $C$ dépendant de $A$ et de $d$. Puis en remplaçant $\alpha$ par $\frac{|x_0-x|^2}{M}$, on voit que pour un $q \in\, ]1, \frac{d}{d-2}[$ si $d > 2$, et $q > 1$ si $d \leq 2$,

$$\int_{\mathbb{R}^d} dx\, (U\nu(x))^q < +\infty$$

Il suffit maintenant d'appliquer l'inégalité de Hölder pour voir que

$$\int_{\mathbb{R}^d} dx\, |\nabla F(x)|^2 U\nu(x) < +\infty$$

lorsque $F$ appartient à $W^{1,2p}$, où $p = \frac{q}{q-1} > 1 \vee \frac{d}{2}$. La proposition 1 permet alors de conclure.

∎

**Corollaire 2** *Si $F$ est dans $H^1_{\text{loc}}$, alors pour quasi tout point de départ $x_0$, le processus $(F(X_t), 0 \leq t \leq S)$ est un processus de Dirichlet sous $\mathbb{P}_{x_0}$ pour tout $S > 0$.*

*Preuve :* En utilisant les estimées inférieures d'Aronson et Nash, et la proposition 3.6 de [1], on voit qu'on peut trouver une suite de compacts $(K_m)_{m\geq 1}$ de $\mathbb{R}^d$, et un ensemble polaire $E_1$ tels que

$$\lim_{m\to+\infty} \mathbb{P}_{x_0}(X_t \in K_m\ \forall t \in [0,S]) = 1 \quad \forall S > 0$$

et

$$\int_{\mathbb{R}^d} |\nabla F_m|^2(x)\, U\delta_{x_0}(x) < +\infty$$

pour tout $x_0 \notin E_1$, tout $m \geq 1$, où $F_m$ désigne la restriction de $F$ à $K_m$ et $\delta_{x_0}$ la mesure de Dirac en $x_0$. D'après la proposition 1, $(Y^m_t := F_m(X_t), 0 \leq t \leq S)$ est un processus de Dirichlet pour tout $m \geq 1$, sous $\mathbb{P}_{x_0}$ pour tout $x_0 \notin E_1$. Sa décomposition s'écrit $Y^m = M^m + N^m$, où $M^m$ est une martingale sous $\mathbb{P}_{x_0}$ pour tout $x_0 \notin E_1$ et $N^m$ est un processus à variations quadratiques nulles sous $\mathbb{P}_{x_0}$ pour tout $x_0 \notin E_1$.

D'après [2], on peut aussi décomposer $Y = F(X)$ en $M + N$ où $M$ est une martingale sous $\mathbb{P}_{x_0}$ pour $x_0 \notin E_2$, où $E_2$ est un ensemble de capacité nulle et $N$ est à variations quadratiques nulles sous toute mesure absolument continue. Soient $T_m = \inf\{t > 0, X_t \notin K_m\}$, $S_n = \sum_{t_i \in D_n}(N_{t_{i+1}} - N_{t_i})^2$ et $S_n^m = \sum_{t_i \in D_n}(N_{t_{i+1}}^m - N_{t_i}^m)^2$. Alors pour $x_0$ en dehors de $E_1 \cup E_2$,

$$\mathbb{P}_{x_0}(S_n > \varepsilon) \leq \mathbb{P}_{x_0}(T_m \leq S) + \mathbb{P}_{x_0}(S_n^m > \varepsilon)$$

D'où

$$\limsup_{n \to +\infty} \mathbb{P}_{x_0}(S_n > \varepsilon) \leq \mathbb{P}_{x_0}(T_m \leq S)$$

Il reste à faire tendre $m$ vers $+\infty$ pour avoir le résultat.

∎

## 2.2  Démonstration de la proposition 1

On ne perd pas en généralité si on prend $S = 1$. Estimons les accroissements du processus $Y$ : pour $T$ et $t$ positifs, on a

$$E_\nu(F(X_{T+t}) - F(X_t))^2 = \int_{\mathbb{R}^d} \nu(dx) \int_{\mathbb{R}^d} dz \int_{\mathbb{R}^d} dy\ p_t(x,y) p_T(y,z) (F(y) - F(z))^2$$

Par ailleurs, pour toute fonction $G$ de classe $C^1$, on a

$$G(z) - G(y) = \int_0^1 \frac{d}{dr} G(y + r(z-y))\, dr = \int_0^1 dr (z-y).\nabla G(y + r(z-y)) \quad (10)$$

Soit $(G_n)_{n \in \mathbb{N}}$ une suite de fonctions de classe $C^1$, telle que

$$\lim_{n \to +\infty} G_n = F \text{ dans } H^1_{\text{loc}}$$

En particulier

$$\lim_{n \to +\infty} G_n(X_t) = F(X_t) \text{ dans } L^2(\mathbb{P}_\nu) \quad \forall t > 0$$

et d'après (10), on peut écrire pour tout $t > 0$ :

$$E_\nu(Y_{T+t} - Y_t)^2 = \lim_{n \to +\infty} E_\nu(G_n(X_{T+t}) - G_n(X_t))^2$$

$$= \lim_{n \to +\infty} \int_{\mathbb{R}^d} \nu(dx) \int_{\mathbb{R}^d} dz \int_{\mathbb{R}^d} dy \, p_t(x,y) p_T(y,z)$$
$$\left( \int_0^1 dr(z-y).\nabla G_n(y + r(z-y)) \right)^2$$

$$\leq \int_{\mathbb{R}^d} \nu(dx) \int_{\mathbb{R}^d} dz \int_{\mathbb{R}^d} dy \, p_t(x,y) p_T(y,z)$$
$$\left( \int_0^1 dr \lim_{n \to +\infty} (z-y).\nabla G_n(y + r(z-y)) \right)^2$$

$$\leq \int_{\mathbb{R}^d} \nu(dx) \int_{\mathbb{R}^d} dz \int_{\mathbb{R}^d} dy \, p_t(x,y) p_T(y,z) \left( \int_0^1 dr(z-y).\nabla F(y + r(z-y)) \right)^2$$

$$\leq \int_{\mathbb{R}^d} \nu(dx) \int_{\mathbb{R}^d} dz \int_{\mathbb{R}^d} dy \, p_t(x,y) p_T(y,z) |y-z|^2 \int_0^1 dr \, |\nabla F(y + r(z-y))|^2$$

en appliquant Cauchy-Schwarz. Le reste du calcul repose sur les estimées d'Aronson et Nash, et des inégalités propres aux probabilités de transition gaussiennes :

$$|y-z|^2 p_T(y,z) \leq |y-z|^2 p_T^M(y,z) \leq 2T \times 2^{d/2} C \times p_{2T}^M(y,z)$$

pour $C = \sup_{x \in \mathbb{R}_+} x \exp -\frac{x}{2}$. Ainsi en renommant la constante,

$$E_\nu(Y_{T+t} - Y_t)^2$$
$$\leq CT \int_{\mathbb{R}^d} \nu(dx) \int_{\mathbb{R}^d} dz \int_{\mathbb{R}^d} dy \, p_t^M(x,y) \, p_{2T}^M(y,z) \int_0^1 dr \, |\nabla F(y + r(z-y))|^2$$
$$\leq CT \int_{\mathbb{R}^d} \nu(dx) \int_{\mathbb{R}^d} dz \int_{\mathbb{R}^d} dy \, p_t^M(x,y) \int_0^1 dr \, p_{2Tr^2}^M(y,z) \, |\nabla F(z)|^2$$

après changement de variable et en notant que $p_{2T}^M(y, \frac{z-y}{r} + y) = r^d p_{2Tr^2}^M(y,z)$

On en déduit

$$E_\nu(Y_{T+t} - Y_t)^2 \leq CT \int_{\mathbb{R}^d} \nu(dx) \int_{\mathbb{R}^d} dz \, |\nabla F(z)|^2 \int_0^1 dr \, p_{t+2Tr^2}^M(x,z) \quad (11)$$

**Lemme 2** *Il existe une constante $C$ telle que pour tout $s \in [t+2T, t+3T]$ on ait :*

$$p_{t+2Tr^2}^M(x,z) \leq C \frac{1}{T} \int_{t+2T}^{t+3T} ds \, p_s^M(x,z) \quad \forall T \leq \frac{t}{2} \quad \forall r \in [0,1]$$

*Preuve :* En effet, pour tout $s \in [t+2T, t+3T]$, pour tout $r \in [0,1]$, on a l'encadrement :
$$1 \leq \frac{s}{t+2Tr^2} \leq \frac{5}{2}$$
Donc $\exp\left(-\frac{|x-z|^2}{M(t+2Tr^2)}\right) \leq \exp\left(-\frac{|x-z|^2}{Ms}\right)$ et $\frac{1}{(t+2Tr^2)^{d/2}} \leq \left(\frac{5}{2}\right)^{d/2} \frac{1}{s^{d/2}}$.
Ce qui donne le résultat pour $C = \left(\frac{5}{2}\right)^{d/2}$

■

En utilisant le lemme 2 et l'inégalité (11), on obtient donc

$$E_\nu(Y_{T+t} - Y_t)^2 \leq C \int_{\mathbb{R}^d} \nu(dx) \int_{\mathbb{R}^d} dz\, |\nabla F(z)|^2 \int_{t+2T}^{t+3T} ds\, p_s^M(x,z) \qquad (12)$$

pour tout $T \leq \frac{t}{2}$.
On applique maintenant l'inégalité (12) à la subdivision $D_n$, pour $T = t_i - t_{i-1}$ et $t = t_{i-1}$. La condition $t_i - t_{i-1} \leq \frac{t_{i-1}}{2}$ est vérifiée pour $i \geq 3$ :

$$\begin{aligned}
E_\nu \sum_{i=3}^{2^n} &\left(Y_{t_i} - Y_{t_{i-1}}\right)^2 \\
&\leq C \int_{\mathbb{R}^d} \nu(dx) \int_{\mathbb{R}^d} dz\, |\nabla F(z)|^2 \sum_{i=3}^{2^n} \int_{t_{i+1}}^{t_{i+2}} ds\, p_s^M(x,z) \\
&\leq C \int_{\mathbb{R}^d} dz\, |\nabla F(z)|^2 \int_0^2 ds \int_{\mathbb{R}^d} \nu(dx) p_s^M(x,z) \\
&\leq C \int_{\mathbb{R}^d} dz\, |\nabla F(z)|^2 U\nu(z)
\end{aligned}$$

De plus d'après (12), $(Y_t)_{t>0}$ est une suite de Cauchy dans $L^2(\mathbb{P}_\nu)$ lorsque $t$ tend vers 0. On peut donc définir $Y_0 := \lim_{t \to 0} Y_t$, dans $L^2(\mathbb{P}_\nu)$. On a alors

$$\lim_{n \to +\infty} E_\nu(Y_{t_2} - Y_0)^2 = 0 \qquad (13)$$

et donc pour $n$ assez grand

$$E_\nu\left(Y_{t_2} - Y_0\right)^2 \leq C \int_{\mathbb{R}^d} dz\, |\nabla F(z)|^2 U\nu(z)$$

*Fin de la preuve :* On conclut comme dans [5] : d'après [2], on peut écrire $Y_t - Y_0 = M_t + N_t$ où $M$ est une martingale sous $\mathbb{P}_{x_0}$ pour tout point $x_0$, car $(X_t, P_t)$ vérifie la condition d'absolue continuité (voir [2], p.208), et $N_t$ est un processus à variations quadratiques nulles sous $\mathbb{P}_\mu$, pour toute mesure de probabilité $\mu$ absolument continue. On doit donc démontrer que $N$ est encore à variations quadratiques nulles sous $\mathbb{P}_\nu$. Comme $P_s\nu$ est absolument continue, par la propriété de Markov, on a :

$$\limsup_{n\to+\infty} \sum_{i=1}^{2^n} \mathbf{1}(t_i \geq s) \left(N_{t_{i+1}} - N_{t_i}\right)^2 = 0 \qquad (14)$$

en probabilité, sous $\mathbb{P}_\nu$, pour tout $0 < s \leq 1$.
De plus d'après (12) et (13), on voit que pour tout $\varepsilon > 0$, on peut choisir $s$ suffisamment petit de manière à obtenir :

$$\limsup_{n\to+\infty} \sum_{i=1}^{2^n} \mathbf{1}(t_{i+1} \leq s) \left(Y_{t_{i+1}} - Y_{t_i}\right)^2 \leq \frac{\varepsilon}{6} \qquad (15)$$

On sait d'autre part que $M$ est une martingale $L^2$, et que

$$\langle M \rangle_t = \int_0^t |\sigma \nabla F|^2(X_s)ds.$$

Ainsi

$$
\begin{aligned}
E_\nu \left(M_{t_{i+1}} - M_{t_i}\right)^2 &= E\left(\nu_{t_i}(X_0) \int_0^{t_{i+1}-t_i} ds \, |\sigma \nabla F|^2(X_s)\right) \\
&= C \int_{\mathbb{R}^d} dx \int_0^{t_{i+1}-t_i} ds \, \nu_{t_i}(x) P_s(|\nabla F|^2)(x) \\
&\leq C \int_{\mathbb{R}^d} dx \int_0^{t_{i+1}-t_i} ds \, \nu_{t_i}(x) P_s^M(|\nabla F|^2)(x) \\
&\leq C \int_{\mathbb{R}^d} dx \int_{t_i}^{t_{i+1}} ds \, \nu_s^M(x) \, |\nabla F|^2(x)
\end{aligned}
$$

D'où l'on déduit :

$$\limsup_{n\to+\infty} \sum_{i=1}^{2^n} \mathbf{1}(t_{i+1} \leq s) \left(M_{t_{i+1}} - M_{t_i}\right)^2 \leq C \int_{\mathbb{R}^d} dx \int_0^s du \, \nu_u^M(x) \, |\nabla F|^2(x) \leq \frac{\varepsilon}{6} \qquad (16)$$

pour $s$ assez petit.
De (15) et (16), on déduit :

$$\limsup_{n\to+\infty} \sum_{i=1}^{2^n} \mathbf{1}(t_{i+1} \leq s) \left(N_{t_{i+1}} - N_{t_i}\right)^2 \leq \frac{2\varepsilon}{3} \qquad (17)$$

Regroupant (14) et (17), il vient

$$\limsup_{n \to +\infty} \mathbb{P}_\nu \left( \sum_{i=1}^{2^n} \left(N_{t_{i+1}} - N_{t_i}\right)^2 \geq a \right)$$
$$\leq \limsup_{n \to +\infty} \mathbb{P}_\nu \left( \sum_{i=1}^{2^n} \mathbf{1}(t_{i+1} \leq s) \left(N_{t_{i+1}} - N_{t_i}\right)^2 \geq \frac{a}{2} \right)$$
$$\leq \limsup_{n \to +\infty} \frac{2}{a} E_\nu \left( \sum_{i=1}^{2^n} \mathbf{1}(t_{i+1} \leq s) \left(N_{t_{i+1}} - N_{t_i}\right)^2 \right)$$
$$\leq \frac{4\varepsilon}{3a}$$

∎

## 3 Formule d'Itô

**Notations 4 :** Par la suite on désignera par $f_k$ (resp. $f_{kl}$) les dérivées faibles $\partial_k F$ (resp. $\partial^2_{kl} F$) d'une fonction $F$, lorsqu'elles existent.

Dans ce paragraphe, on démontre une formule d'Itô sous $\mathbb{P}_\nu$ lorsque $F$ et $\nu$ sont liées par

$$\int dz |\nabla F(z)|^2 U\nu(z) < +\infty \tag{18}$$

$$\int_{\mathbb{R}^d} dz \sum_{k,l=1}^d f_{kl}^2(z) U\nu(z) < +\infty \tag{19}$$

Pour cela on démontre d'abord l'existence des covariations quadratiques $[f(X), X]$ lorsque $f$ vérifie

$$\int dz |\nabla f(z)|^2 U\nu(z) < +\infty \tag{20}$$

en montrant la convergence en probabilité des sommes

$$\sum_{t_i \in D_n} \left(f(X_{t_{i+1}}) - f(X_{t_i})\right) \left(X_{t_{i+1}}^k - X_{t_i}^k\right)$$

On montre de plus que ces sommes sont dominées dans $L^1(\mathbb{P}_\nu)$ par la norme définie par (20). En seconde partie, on définit l'intégrale progressive $\int_0^S f(X_s)\, dX_s$, pour toute fonction $f$ satisfaisant (20), comme limite des sommes de Riemann $\sum_{t_i \in D_n} f(X_{t_i})(X_{t_{i+1}} - X_{t_i})$ en probabilité, sous $\mathbb{P}_\nu$.

Le calcul est fait lorsque $f$ est une dérivée faible $f = \partial_k F$ où $F$ vérifie (18) et (19)

L'idée est d'estimer dans $L^1(\mathbb{P}_\nu)$ les sommes

$$\sum_{t_i \in D_n} \left( F(X_{t_{i+1}}) - F(X_{t_i}) - \nabla F(X_{t_i}).(X_{t_{i+1}} - X_{t_i}) \right)$$

grâce encore à une formule de Taylor, puis de déduire leur convergence du cas où $F$ est une fonction de classe $C^2$.

## 3.1 Covariations quadratiques

**Proposition 2** *Soit $f$ une fonction satisfaisant la condition (20). Alors pour tout $k = 1, \ldots, d$*

$$[f(X), X^k]_S := \lim_{n \to +\infty} \sum_{t_i \in D_n} (f(X_{t_{i+1}}) - f(X_{t_i}))(X^k_{t_{i+1}} - X^k_{t_i})$$

*existe en probabilité sous $\mathbb{P}_\nu$.*
*De plus on a pour $n$ assez grand,*

$$E_\nu \left( \sum_{t_i \in D_n} \left| \left(f(X_{t_{i+1}}) - f(X_{t_i})\right) \left(X^k_{t_{i+1}} - X^k_{t_i}\right) \right| \right)$$
$$\leq C \left( \int dz |\nabla f(z)|^2 U\nu(z) \right)^{1/2}$$

*Preuve :* 1) D'après la première partie, on sait qu'on peut trouver une version quasi-continue de $f$ telle que $f(X_0)$ soit "bien" défini. Comme $f$ vérifie (18), $f(X)$ est un processus de Dirichlet d'après la proposition 1 et s'écrit $f(X_t) = M^f_t + N^f_t$. $M^f$ est une martingale sous $\mathbb{P}_\nu$ ; $N^f$ est un processus à variations quadratiques nulles sous $\mathbb{P}_\nu$.
De même, soit $M^k + N^k$ la décomposition de $X^k$ en tant que processus de Dirichlet local (l'application $e_k : x \mapsto x_k$ appartient à $H^1_{\text{loc}}$). $M^k$ est une martingale $L^2$ car $E_\nu \langle M^k \rangle_t = \int \nu(dx) \int_0^t ds \, P_s(|\sigma \nabla e_k|^2)(x) < +\infty$. On a alors

$$\sum_{t_i \in D_n} (f(X_{t_{i+1}}) - f(X_{t_i}))(X^k_{t_{i+1}} - X^k_{t_i})$$
$$= \sum_{t_i \in D_n} (M^f_{t_{i+1}} - M^f_{t_i})(M^k_{t_{i+1}} - M^k_{t_i}) + \sum_{t_i \in D_n} (M^f_{t_{i+1}} - M^f_{t_i})(N^k_{t_{i+1}} - N^k_{t_i})$$
$$+ \sum_{t_i \in D_n} (N^f_{t_{i+1}} - N^f_{t_i})(M^k_{t_{i+1}} - M^k_{t_i}) + \sum_{t_i \in D_n} (N^f_{t_{i+1}} - N^f_{t_i})(N^k_{t_{i+1}} - N^k_{t_i})$$

Chacun des trois derniers termes du membre de droite de l'égalité a une limite nulle, dans $L^1(\mathbb{P}_\nu)$, lorsque $n$ croît indéfiniment en vertu de l'inégalité de Cauchy-Schwarz. Le premier terme tend vers le crochet $\langle M^f, M^k \rangle_S$, dans $L^1(\mathbb{P}_\nu)$.

2) On montre maintenant l'estimation de ces sommes : d'une part, on a

$$E_\nu \left( \sum_{\substack{t_i \in D_n \\ i \geq 2}} (X^k_{t_{i+1}} - X^k_{t_i})^2 \right) = \sum_{\substack{t_i \in D_n \\ i \geq 2}} \int_{\mathbb{R}^d} dx\, \nu_{t_i}(x) \int_{\mathbb{R}^d} dy\, p_{t_{i+1}-t_i}(x,y)(x_k - y_k)^2$$

$$\leq C \sum_{\substack{t_i \in D_n \\ i \geq 2}} \int_{\mathbb{R}^d} dx\, \nu_{t_i}(x) \int_{\mathbb{R}^d} dy\, p^M_{t_{i+1}-t_i}(x,y)(x_k - y_k)^2$$

avec $\nu_{t_i}(x) \leq C \nu^M_{t_i}(x)$. Puis, si $t_i \leq s \leq t_{i+1}$,

$$p^M_{t_i}(x,y) = \frac{M}{t_i^{d/2}} \exp -\left( \frac{|x-y|^2}{Mt_i} \right)$$

$$\leq \left(\frac{s}{t_i}\right)^{d/2} \frac{M}{s^{d/2}} \exp -\left( \frac{|x-y|^2}{Ms} \right)$$

$$\leq \sup_{i>0} \left(\frac{t_{i+1}}{t_i}\right)^{d/2} \frac{M}{s^{d/2}} \exp -\left( \frac{|x-y|^2}{Ms} \right)$$

$$\leq 2^{d/2} p^M_s(x,y)$$

D'où $\nu_{t_i}(x) \leq C \frac{1}{t_{i+1}-t_i} \int_{t_i}^{t_{i+1}} \nu^M_s(x)\, ds$. De plus,

$$\int_{\mathbb{R}^d} dy\, p^M_{t_{i+1}-t_i}(x,y)(x_k - y_k)^2 = C(t_{i+1} - t_i)$$

D'où

$$E_\nu \left( \sum_{\substack{t_i \in D_n \\ i \geq 2}} (X^k_{t_{i+1}} - X^k_{t_i})^2 \right) \leq C \sum_{\substack{t_i \in D_n \\ i \geq 2}} \int_{\mathbb{R}^d} dx \int_{t_i}^{t_{i+1}} ds\, \nu^M_s(x)$$

$$\leq C \int_{\mathbb{R}^d} dx \int_0^S ds\, \nu^M_s(x) \leq C$$

car $\int_{\mathbb{R}^d} dx\, \nu^M_s(x) = C$ où $C$ est une autre constante qui ne dépend que de $M$ et de la dimension $d$.

Et comme $\lim_{n \to +\infty} E_\nu(X^k_{t_2} - X^k_0)^2 = 0$, ceci reste vrai si on somme sur toute la subdivison $D_n$.

D'autre part, on a d'après l'estimation (12)

$$E_\nu \left( \sum_{t_i \in D_n} (f(X_{t_{i+1}}) - f(X_{t_i}))^2 \right)$$

$$\leq C \int dz |\nabla f(z)|^2 \int_0^{2S} ds\, \nu^M_s(z) + E_\nu(f(X_{t_2}) - f(X_0))^2$$

Comme de plus $\lim_{n\to+\infty} E_\nu(f(X_{t_2}) - f(X_0))^2 = 0$, on peut toujours supposer que pour $n$ assez grand on a aussi

$$E_\nu\left((f(X_{t_2}) - f(X_0))^2\right) \leq C \int dz |\nabla f(z)|^2 \int_0^{2S} ds\, \nu_s^M(z)$$

D'où pour $n$ assez grand, en utilisant Cauchy-Schwarz, on déduit :

$$E_\nu\left(\sum_{t_i \in D_n} \left|(f(X_{t_{i+1}}) - f(X_{t_i}))(X_{t_{i+1}}^k - X_{t_i}^k)\right|\right)$$

$$\leq E_\nu\left(\sum_{t_i \in D_n}(f(X_{t_{i+1}}) - f(X_{t_i}))^2\right)^{\frac{1}{2}} \left(\sum_{t_i \in D_n}(X_{t_{i+1}}^k - X_{t_i}^k)^2\right)^{\frac{1}{2}}$$

$$\leq \left(E_\nu \sum_{t_i \in D_n}(f(X_{t_{i+1}}) - f(X_{t_i}))^2\right)^{\frac{1}{2}} \left(E_\nu \sum_{t_i \in D_n}(X_{t_{i+1}}^k - X_{t_i}^k)^2\right)^{\frac{1}{2}}$$

$$\leq C \left(\int dz |\nabla f(z)|^2 \int_0^{2S} ds\, \nu_s^M(z)\right)^{\frac{1}{2}}$$

$$\leq C \left(\int dz |\nabla f(z)|^2 U\nu(z)\right)^{\frac{1}{2}}$$

∎

## 3.2 Intégrale progressive

**Proposition 3** *On se donne une fonction $F$ vérifiant les conditions (18) et (19) alors on a l'estimation suivante pour $n$ assez grand :*

$$E_\nu\left(\sum_{t_i \in D_n} |F(X_{t_{i+1}}) - F(X_{t_i}) - \nabla F(X_{t_i}).(X_{t_{i+1}} - X_{t_i})|\right)$$

$$\leq C \sum_{k,l=1}^d \left(\int_{\mathbb{R}^d} dz\, f_{kl}^2(z) U\nu(z)\right)^{\frac{1}{2}}$$

D'après la remarque faite à la proposition 1, $F(X_0)$ et $\nabla F(X_0)$ sont bien définis.

*Preuve :* Pour $T \geq 0$ et $t > 0$, on a :

$$E_\nu |F(X_{T+t}) - F(X_t) - \nabla F(X_t).(X_{T+t} - X_t)|$$
$$= \iint dx\, dy\, \nu_t(x) p_T(x,y) |F(y) - F(x) - \nabla F(x).(y-x)|$$
$$\leq C \iint dx\, dy\, \nu_t^M(x) p_T^M(x,y) |F(y) - F(x) - \nabla F(x).(y-x)|$$

Pour une fonction $G$ de classe $C^2$, on a

$$G(y) - G(x) - \nabla G(x).(y-x) = \int_0^1 (y-x).\left(\nabla G(x+r(y-x)) - \nabla G(x)\right)$$

$$= \sum_{k,l=1}^d \int_0^1 r\,dr \int_0^1 ds\, (y-x)_k\,(y-x)_l\,g_{kl}(x+sr(y-x))$$

Les conditions (18) et (19) impliquent que $F$ appartient à $H^2_{\text{loc}}$. Soit $(G_n)_{n\in\mathbb{N}}$ une suite de fonctions $C^2$ qui convergent vers $F$ dans $H^2_{\text{loc}}$. En passant à la limite dans
$E_\nu \left|G_n(X_{T+t}) - G_n(X_t) - \nabla G_n(X_t).(X_{T+t} - X_t)\right| \leq$
$C \sum_{k,l=1}^d \iint dx\,dy\,\nu_t^M(x) p_T^M(x,y) \int_0^1 r\,dr \int_0^1 ds\, \left|(y-x)_k\,(y-x)_l\,g_{kl}^{(n)}(x+sr(y-x))\right|$
où $g_{kl}^{(n)} = \partial_{kl}^2 G_n$

On obtient
$E_\nu \left|F(X_{T+t}) - F(X_t) - \nabla F(X_t).(X_{T+t} - X_t)\right| \leq$
$\sum_{k,l=1}^d \iint dx\,dy\,\nu_t^M(x) p_T^M(x,y) \int_0^1 r\,dr \int_0^1 ds\, |(y-x)_k\,(y-x)_l\,f_{kl}(x+sr(y-x))|$

Dans l'intégrale en $y$, on note que

$$|x-y|_k\,|x-y|_l\,p_T^M(x,y) \leq CT\,p_{2T}^M(x,y)$$

pour une constante $C$, et on effectue le changement de variable $z = x + sr(y-x)$.
Il vient :

$$\int dy\,p_T^M(x,y)|x-y|_k\,|x-y|_l\,|f_{kl}(x+sr(y-x))|$$

$$\leq CT \int dz\,\frac{1}{(sr)^d}\,p_{2T}^M\left(x, \frac{z-x}{sr}+x\right)|f_{kl}(z)|$$

$$\leq CT \int dz\,p_{2T(sr)^2}^M(x,z)\,|f_{kl}(z)|$$

D'où

$$E_\nu \left|F(X_{T+t}) - F(X_t) - \nabla F(X_t).(X_{T+t} - X_t)\right|$$

$$\leq CT \sum_{k,l=1}^d \iint dx\,dy \int_0^1 r\,dr \int_0^1 ds\,\nu_t^M(x) p_{2T(sr)^2}^M(x,y)\,|f_{kl}(y)|$$

$$\leq CT \sum_{k,l=1}^d \iint \nu(dx)\,dy \int_0^1 du \int_0^1 dv\,p_{2Tu^2+t}^M(x,y)\,|f_{kl}(y)|$$

$$\leq C \sum_{k,l=1}^d \int \nu(dx) \int dy \int_{t+2T}^{t+3T} ds\,p_s^M(x,y)\,|f_{kl}(y)| \qquad (21)$$

La dernière inégalité résulte du lemme 2, pour $T \leq \frac{t}{2}$.
Appliquons (21) à la subdivision $D_n$, pour $T = t_{i+1} - t_i$ et $t = t_i$, $i \geq 2$ :

$$E_\nu \left( \sum_{\substack{t_i \in D_n \\ i \geq 2}} |F(X_{t_{i+1}}) - F(X_{t_i}) - \nabla F(X_{t_i}).(X_{t_{i+1}} - X_{t_i})| \right)$$

$$\leq C \sum_{k,l=1}^d \int dy \, |f_{kl}|(y) \int_0^{2S} ds \, \nu_s^M(y)$$

$$\leq C \sum_{k,l=1}^d \left( \int dy \, |f_{kl}|^2(y) \int_0^{2S} ds \, \nu_s^M(y) \right)^{\frac{1}{2}} \left( \int dy \int_0^{2S} ds \, \nu_s^M(y) \right)^{\frac{1}{2}}$$

$$\leq C \sum_{k,l=1}^d \left( \int dy \, |f_{kl}|^2(y) \int_0^{2S} ds \, \nu_s^M(y) \right)^{\frac{1}{2}}$$

Cette inégalité est encore vraie si on somme sur les indices $i = 0, 1$, pour $n$ assez grand, car de (21) on déduit aussi que

$$\lim_{n \to +\infty} E_\nu |F(X_{t_2}) - F(X_0) - \nabla F(X_0).(X_{t_2} - X_0)| = 0. \qquad \blacksquare$$

On peut maintenant en déduire le théorème suivant :

**Théorème 1** *Soit $F$ une fonction de $L^2(dx)$ qui vérifie (18) et (19). Alors,*

$$\int_0^S \nabla F(X_s) dX_s := \lim_{n \to +\infty} \sum_{t_i \in D_n} \nabla F(X_{t_i}).(X_{t_{i+1}} - X_{t_i}) \qquad (22)$$

*existe en probabilité sous $\mathbb{P}_\nu$. De plus on a la formule d'Itô*

$$F(X_S) - F(X_0) = \int_0^S \nabla F(X_s) dX_s + \frac{1}{2} \sum_{k=1}^d [f_k(X), X^k]_S \quad \mathbb{P}_\nu - p.s. \qquad (23)$$

*Preuve :*
1) Comme $F$ vérifie (18) et (19), elle appartient à $H^2_{\text{loc}}$. On commence par montrer que si $F$ est de classe $C_c^2$, alors la limite (22) existe sous $\mathbb{P}_\nu$. Ensuite on approchera $F$ dans $H^2_{\text{loc}}$ par une suite $(G_p)_{p \in \mathbb{N}}$ de fonctions $C_c^2$ sur $\mathbb{R}^d$, et on passera à la limite grâce aux estimations faites précédemment.
Soit $G$ une fonction de classe $C_c^2$ sur $\mathbb{R}^d$. Pour tout $x$ et tout $y$ de $\mathbb{R}^d$, on a :

$$|G(x+y) - G(x) - \frac{1}{2} \left( \nabla G(x).y + \nabla G(x+y).y \right)| \leq \sup_{z \in [x, x+y]} |D^2 G(z)| \, |y|^2$$

En appliquant ceci à $x+y = X_{t_{i+1}}$ et $x = X_{t_i}$, et en sommant sur les $t_i$, on obtient

$$E_\nu \left( \sum_{t_i \in D_n} |G(X_{t_{i+1}}) - G(X_{t_i}) - \frac{1}{2}(\nabla G(X_{t_i}) + \nabla G(X_{t_{i+1}})).(X_{t_{i+1}} - X_{t_i})| \right)$$

$$\leq E_\nu \left( \sum_{t_i \in D_n} \sup_{z \in [X_{t_i}, X_{t_{i+1}}]} |D^2 G(z)| \, |X_{t_{i+1}} - X_{t_i}|^2 \right)$$

$$\leq E_\nu \left( \sup_{t_i \in D_n} \sup_{z \in [X_{t_i}, X_{t_{i+1}}]} |D^2 G(z)| \sum_{t_i \in D_n} |X_{t_{i+1}} - X_{t_i}|^2 \right)$$

D'où

$$\liminf_{n \to +\infty} E_\nu \left( \sum_{t_i \in D_n} |G(X_{t_{i+1}}) - G(X_{t_i}) - \frac{1}{2}(\nabla G(X_{t_i}) + \nabla G(X_{t_{i+1}})).(X_{t_{i+1}} - X_{t_i})| \right)$$

$$\leq E_\nu \left( \liminf_{n \to +\infty} \sup_{t_i \in D_n} \sup_{z \in [X_{t_i}, X_{t_{i+1}}]} |D^2 G(z)| \sum_{t_i \in D_n} |X_{t_{i+1}} - X_{t_i}|^2 \right)$$

$$= 0$$

par continuité des trajectoires $s \mapsto G(X_s)$, et parce que $\sum_{t_i \in D_n} |X_{t_{i+1}} - X_{t_i}|^2$ est borné dans $L^2(\mathbb{P}_\nu)$.

Comme

$$\lim_{n \to +\infty} \frac{1}{2} E_\nu \left( \sum_{t_i \in D_n} (\nabla G(X_{t_{i+1}}) - \nabla G(X_{t_i})).(X_{t_{i+1}} - X_{t_i}) \right)$$

existe d'après la proposition 2, on en déduit que

$$\lim_{n \to +\infty} E_\nu \left( \sum_{t_i \in D_n} G(X_{t_{i+1}}) - G(X_{t_i}) - \nabla G(X_{t_i}).(X_{t_{i+1}} - X_{t_i}) \right)$$

existe aussi et que ces deux limites sont égales. En particulier

$$\lim_{n \to +\infty} \left( \sum_{t_i \in D_n} \nabla G(X_{t_i}).(X_{t_{i+1}} - X_{t_i}) \right)$$

existe en probabilité sous $\mathbb{P}_\nu$ et cette limite, notée $\int_0^S G(X_s) dX_s$ est égale à

$$\int_0^S G(X_s) dX_s = G(X_S) - G(X_0) - \frac{1}{2} \sum_{k=1}^d [g_k(X), X^k]_S \quad \mathbb{P}_\nu - \text{p.s.}$$

2) Soit $(G_p)$ une suite de fonctions $C_c^2$ qui converge vers $F$ dans $H_{\text{loc}}^2$. Alors on a d'une part

$$\lim_{p\to+\infty} \int_{\mathbb{R}^d} dz\, |\nabla F(z) - \nabla G_p(z)|^2\, U\nu(z) = 0 \qquad (24)$$

$$\lim_{p\to+\infty} \int_{\mathbb{R}^d} dz\, \sum_{k,l=1}^d (f_{kl}(z) - g_{kl}^{(p)}(z))^2\, U\nu(z) = 0 \qquad (25)$$

et d'autre part

$$\lim_{n\to+\infty} E_\nu\left(\sum_{t_i \in D_n} G_p(X_{t_{i+1}}) - G_p(X_{t_i}) - \nabla G_p(X_{t_i}).(X_{t_{i+1}} - X_{t_i})\right) =$$

$$\lim_{n\to+\infty} \frac{1}{2} E_\nu\left(\sum_{t_i \in D_n} (\nabla G_p(X_{t_{i+1}}) - \nabla G_p(X_{t_i})).(X_{t_{i+1}} - X_{t_i})\right) \qquad (26)$$

Grâce aux estimées des propositions 2 et 3, on peut passer à la limite en $p$ dans (26), et les égalités (24) et (25) prouvent qu'à la limite on a

$$\lim_{n\to+\infty} \left(\sum_{t_i \in D_n} F(X_{t_{i+1}}) - F(X_{t_i}) - \nabla F(X_{t_i}).(X_{t_{i+1}} - X_{t_i})\right) =$$

$$\lim_{n\to+\infty} \frac{1}{2}\left(\sum_{t_i \in D_n} (\nabla F(X_{t_{i+1}}) - \nabla F(X_{t_i})).(X_{t_{i+1}} - X_{t_i})\right) \quad \mathbb{P}_\nu - \text{p.s.}$$

En particulier

$$\int_0^S \nabla F(X_s) dX_s := \lim_{n\to+\infty} \sum_{t_i \in D_n} \nabla F(X_{t_i}).(X_{t_{i+1}} - X_{t_i})$$

existe en probabilité sous $\mathbb{P}_\nu$ et on a :

$$F(X_S) - F(X_0) - \int_0^S \nabla F(X_s) dX_s = \frac{1}{2}\sum_{k=1}^d [f_k(X), X^k]_S \quad \mathbb{P}_\nu - \text{p.s.}$$

ce qui termine la preuve.

∎

**Corollaire 3** *Soit $F$ une fonction de $H_{\text{loc}}^2$. Pour quasi tout point $x_0$,*

$$\int_0^S \nabla F(X_s) dX_s := \lim_{n\to+\infty} \sum_{t_i \in D_n} \nabla F(X_{t_i}).(X_{t_{i+1}} - X_{t_i})$$

*existe en probabilité sous $\mathbb{P}_{x_0}$. De plus on a la formule d'Itô*

$$F(X_t) - F(X_0) = \int_0^S \nabla F(X_s) dX_s + \frac{1}{2}\sum_{k=1}^d [f_k(X), X^k]_S \quad \mathbb{P}_{x_0} - p.s.$$

*Preuve :* En utilisant les arguments du corollaire 2, on peut toujours suppposer que $F$ est à support compact et qu'elle vérifie (18) et (19) pour $\nu = \delta_{x_0}$, la mesure de Dirac en $x_0$, pour quasi tout $x_0$. Alors d'après le théorème 1, la limite

$$\lim_{n \to +\infty} \sum_{t_i \in D_n} \nabla F(X_{t_i}).(X_{t_{i+1}} - X_{t_i})$$

existe en probabilité sous $\mathbb{P}_{x_0}$ pour quasi tout $x_0$, et on a la formule d'Itô $\mathbb{P}_{x_0}$ – p.s..

■

**Corollaire 4** *Si $F \in W^{2,p}$ pour un $p > d \vee 2$*

$$\int_0^S \nabla F(X_s) dX_s = \lim_{n \to +\infty} \sum_{t_i \in D_n} \nabla F(X_{t_i}).(X_{t_{i+1}} - X_{t_i})$$

*et*

$$\sum_{k=1}^d [f_k(X), X^k]_t = \lim_{n \to +\infty} \sum_{k=1}^d \sum_{t_i \in D_n} (f_k(X_{t_{i+1}}) - f_k(X_{t_i}))(X_{t_{i+1}}^k - X_{t_i}^k)$$

*existent en probabilité sous $\mathbb{P}_{x_0}$ pour tout $x_0$. Et on a la formule d'Itô (23).*

*Preuve :* Si $F \in W^{2,p}$ pour un $p > d \vee 2$ choisi comme dans le corollaire 1, d'après l'inégalité de Hölder, les conditions (18) et (19) sont vérifiées pour $\nu = \delta_{x_0}$, pour tout $x_0 \in \mathbb{R}^d$. Ceci donne l'existence de l'intégrale progressive $\displaystyle\int_0^S \nabla F(X_s) dX_s$ comme limite en probabilité des sommes de Riemann $\displaystyle\sum_{t_i \in D_n} \nabla F(X_{t_i}).(X_{t_{i+1}} - X_{t_i})$ sous $\mathbb{P}_{x_0}$ pour tout $x_0$, et l'existence des covariations quadratiques $[f_k(X), X^k]_S$ sous $\mathbb{P}_{x_0}$ pour tout $k = 1, \ldots, d$. Ainsi la formule d'Itô (23) a lieu $\mathbb{P}_{x_0}$ – p.s. pour tout $x_0$.

■

# Références

[1] H. FÖLLMER, P. PROTTER(1999) On Itô's formula for multidimensionnal Brownian motion, *Probability theory and related fields No 116* pp. 1-20

[2] M. FUKUSHIMA, Y. OSHIMA, M. TAKEDA(1994) Dirichlet forms and symmetric Markov processes, *De Gruyter studies in math. 19*

[3] R. HÖHNLE, K.-TH. STURM(1993) A multidimensionnal analogue to the 0-1 law of Engelbert and Schmidt, *Stochastics and stochastics reports No 44* pp. 27-41

[4] T.J. LYONS, T.S. ZHANG(1994) Decomposition of Dirichlet processes and its application, *The annals of probability Vol.22, No 1* pp.494-524

[5] P. MATHIEU(2000), Dirichlet processes associated to diffusions, *Stochastics and stochastics reports No 71* pp.165-176

[6] A. ROZKOSZ(1996) Stochastic representation of diffusions corresponding to divergence form operators, *Stochastic processes and their Applications No 63* pp. 11-33

[7] D. W. STROOCK(1988) Diffusion semigroups corresponding to uniformly elliptic divergence form operator, *Séminaire de Probabilités XXII, Lecture Notes in Mathematics No 1321* pp. 316-347